\documentclass[12pt]{amsart}
\usepackage{amsfonts}
\usepackage{amssymb}
\usepackage{amsmath}
\usepackage{graphicx}
\usepackage{hyperref}

\def\be{\begin{equation}}
\def\ee{\end{equation}}
\def\bea{\begin{eqnarray*}}
\def\eea{\end{eqnarray*}}

\def\Ric{r}

\DeclareMathOperator{\Aut}{Aut}

\DeclareMathOperator{\Rm}{Rm}

\newtheorem{thm}{Theorem}

\def\RR{{\mathbb R}}
\def\CP{{\mathbb C \mathbb P}}

\begin{document}
\author{Caner Koca}

\address{Stony Brook University, Department of Mathematics, Stony Brook, NY 11794, USA}
\email{caner@math.sunysb.edu}

\title{Einstein Hermitian Metrics of Positive Sectional Curvature}

\begin{abstract}
It is shown that, up to scaling and isometry, the only complete 4-manifold with an Einstein metric of positive sectional curvature which is also Hermitian with respect to some complex structure, is the complex projective plane $\CP_2$, with its Fubini-Study metric.
\end{abstract}
\maketitle

\section{Introduction}

Let $M$ be a smooth $n$-manifold. A Riemannian metric $g$ on $M$ is called \emph{Einstein} if the Ricci tensor is a constant multiple of the metric tensor, i.e.
\[
\Ric=\lambda g
\]
for some constant $\lambda\in\RR$, called the \emph{Einstein constant} \cite{Bes:1987}. If $\lambda>0$ and $g$ is complete, $M$ is compact by Myer's Theorem. 

Since the Ricci tensor is, by definition, the pointwise average of all sectional curvatures, the positivity of $\lambda$ is guaranteed for any Einstein manifold of positive sectional curvature. This paper will investigate the existence of positive sectional curvature Einstein metrics on 4-manifolds. Examples of such manifolds are $4$-sphere $S^4$ with its standard round metric (which has all sectional curvatures $K\equiv 1$), and the complex projective plane $\CP_2$ with the Fubini-Study metric $g_{FS}$ (which has $1\leq K \leq 4$ everywhere). Notice that these examples are of strictly \emph{positive} sectional curvature. The product metric on $S^2\times S^2$ is Einstein with $\lambda>0$, too, its sectional curvatures, however, are \emph{non-negative}. They are actually $0$ for the tangent planes which have 1-dimensional projections on each factors. In fact, the famous Hopf Conjecture asks whether or not there are any metrics on $S^2\times S^2$ of positive sectional curvature.

One fruitful source of Einstein metrics is provided by the K\"ahler geometry. For a compact complex surface $M$,
\begin{enumerate}
 \renewcommand{\labelenumi}{\textbf{\alph{enumi}.}}
 \item there is a unique K\"ahler-Einstein metric with $\lambda<0$ if $c_1(M)<0$ (see \cite{Y2:1978}, \cite{Aub:1982}),
 \item there is a unique K\"ahler-Einstein metric with $\lambda=0$ (i.e. Ricci-flat) in each K\"ahler class if  $c_1(M)=0$ (see \cite{Y1:1977}).
 \item In $\lambda>0$ case, Tian \cite{T:1990} showed that $M$ admits a K\"ahler-Einstein metric with $\lambda>0$ iff $M$ has $c_1(M)>0$ and its automorphism group $\Aut(M)$ is a reductive Lie group. The diffeomorphism types of such complex surfaces are $\CP_2$, $\CP_1\times \CP_1$ and $\CP_2\sharp k\overline{\CP}_2$ with $k=3,4,\dots,8$.
\end{enumerate}

Among these K\"ahler-Einstein metrics, only the one on $\CP_2$ has \emph{positive} sectional curvature.  This follows, for example, by Frankel's Theorem \cite{Fra:1961}, which says that any compact complex surface $M$ with a K\"ahler metric of positive sectional curvature must be $\CP_2$.\footnote[2]{Gursky and LeBrun have proved a stronger result which applies to more general Einstein metrics. See \cite{GL:1999}, Theorem B.} This result, which was proved using a result of Andreotti \cite{And:1957}, is a special case of the Frankel conjecture \cite{Fra:1961}, proved by Siu and Yau \cite{SY:1980}, which asserts that the generalization of the statement is true for all complex dimensions $n$.

If we relax the K\"ahler condition on the Einstein metric $g$, and merely assume that $g$ is \emph{Hermitian}, that is $g(J\cdot,J\cdot)=g(\cdot,\cdot)$ for a complex structure $J$ on the manifold $M$, quite remarkably, only two additional exceptional metrics appear:

\begin{thm}[LeBrun \cite{LeB:2011}] \label{LeBrun:classification}Let $(M^4,J)$ be a compact complex surface. If $g$ is Einstein and Hermitian, then only one of the following holds:
\begin{enumerate}
 \renewcommand{\labelenumi}{(\arabic{enumi})}
\item $g$ is K\"ahler-Einstein with $\lambda >0$.
\item $(M,J)$ is biholomorphic to $\CP_2\sharp\overline{\CP}_2$ and $g$ is the Page metric $g_{\mathrm{Page}}$ (up to rescaling and isometry).
\item $(M,J)$ is biholomorphic to $\CP_2\sharp 2\overline{\CP}_2$ and $g$ is the Chen-LeBrun-Weber metric $g_{\mathrm{CLW}}$ (up to rescaling and isometry).
\end{enumerate}
\end{thm}

If we futhermore assume that $g$ has positive sectional curvature, then the first case of the above theorem is possible only when $M$ is $\CP_2$, by Frankel's theorem. On the other hand, by a theorem of Berger \cite{Ber:1963}, the Fubini-Study metric is the only K\"ahler-Einstein metric with positive holomorphic bisectional\footnote[1]{For a Hermitian metric, the \emph{holomorphic bisectional curvature} of a pair of unit tangent vectors $X,Y$ is defined as $H(X,Y)=\Rm(X,JX,Y,JY)$. Moreover, if the metric is \emph{K\"ahler}, then $H(X,Y)=K(X,Y)+K(X,JY)$, where $K(\cdot,\cdot)$ stands for the sectional curvature. Therefore, positivity of sectional curvature implies that of bisectional curvature.} curvature on $\CP_2$ (up to rescaling and isometry). Therefore, it remains to prove that $g_{\mathrm{Page}}$ and $g_{\mathrm{CLW}}$ are not of positive sectional curvature either. This will conclude the proof of our main theorem:

\begin{thm}\label{maintheorem}
Let $M$ be a smooth $4$-manifold, and let $g$ be a complete Einstein metric of positive sectional curvature. If $g$ is Hermitian with respect to some complex structure $J$ on $M$, then $(M,J)$ is biholomorphic to $\CP_2$, and $g$ is the Fubini-Study metric (up to rescaling and isometry).
\end{thm}

The key fact we use in the proof of Theorem \ref{maintheorem} is  a result due to Frankel \cite{Fra:1961} which says that totally geodesic submanifolds of complementary dimensions on positively curved manifolds necessarily intersect. The strategy of our proof is to obtain a pair of disjoint totally geodesic surfaces for each of the exceptional metrics of Theorem \ref{LeBrun:classification}.

\section{Proof of the Main Theorem}
Let $(M^4,J,g)$ be an Einstein Hermitian manifold of positive sectional curvature. If the metric is \emph{K\"ahler}, by our discussion in the Introduction, $(M,J,g)$ is $\CP_2$ with the Fubini-Study metric (up to rescaling and isometry).

If $g$ is \emph{not} K\"ahler, it follows from LeBrun's Theorem that, it is either $\CP_2\sharp\overline{\CP}_2$ with $g_{\mathrm{Page}}$ or $\CP_2\sharp 2\overline{\CP}_2$ with $g_{\mathrm{CLW}}$. Both of these metrics are conformally equivalent to \emph{extremal K\"ahler metrics} \cite{LeB:1995}. On the other hand, Calabi showed \cite{C2:1985} that the identity component of the isometry group of an extremal K\"ahler metric necessarily lies in the identity component of the group of biholomorphisms as a maximal compact subgroup. The crucial point is that for these two metrics in question, the isometry group contains a torus \footnote[3]{For the Page metric, the isometry group is indeed $U(2)$; see\cite{BB:1982}.} (see \cite{BB:1982}, \cite{CLW:2008}). It is a well-known fact in toric geometry that the images of the moment maps (see \cite{Ful:1993}, \cite{CLW:2008}) of these torus actions are given as follows:\\
\begin{center}
\begin{tabular}{ccc}
\includegraphics[scale=1]{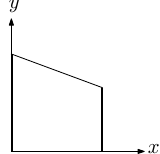} &\quad\quad\quad &\includegraphics[scale=1]{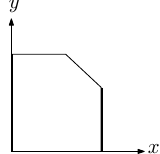} \\
$\CP_2\sharp\overline{\CP}_2$ & &$\CP_2\sharp 2\overline{\CP}_2$
\end{tabular}
\end{center}
where $x$ and $y$ are the Hamiltonians of the torus action. The two vertical edges of each of these moment map profiles correspond to a pair of disjoint 2-spheres in the respective manifolds\footnote[4]{In the 1-point blowup case, these 2-spheres are precisely the 0- and $\infty$-sections of the Hirzebruch Surface $\mathbb F_1$.}, which are fixed by the circle action generated by the Hamiltonian vector field corresponding to the Hamiltonian $x$. Note that this Hamiltonian action preserves the K\"ahler form and the complex structure. As a consequence, we obtain an $S^1$-family of isometries of the K\"ahler metric and the conformally equivalent Einstein metric fixing the two disjoint 2-spheres. Because the common fixed point set of a family of isometries is totally geodesic (see \cite{K:1958}), these submanifolds are moreover totally geodesic.

On the other hand, using a variational argument due to Synge \cite{Synge:1936}, Frankel proved the following theorem about totally geodesic submanifolds of positively curved spaces:
\begin{thm}[Frankel \cite{Fra:1961}]
Let $M$ be a smooth $n$-manifold, and let $g$ be a complete Riemannian metric of positive sectional curvature. If $X$ and $Y$ are two compact totally geodesic submanifolds of dimensions $d_1$ and $d_2$ such that $d_1+d_2\geq n$, then $X$ and $Y$ intersect.
\end{thm}
Since we have found two disjoint totally geodesic submanifolds of dimension 2, neither $g_{\mathrm{Page}}$ nor $g_{\mathrm{CLW}}$ can have positive sectional curvature. This completes the proof of Theorem \ref{maintheorem}.

\bibliography{page.pams.version}
\bibliographystyle{abbrv}

\noindent\textbf{Acknowledgements.} I would like to thank my advisor Claude LeBrun for suggesting the problem and for his help, guidance and encouragement. Thanks for the referee for many useful comments and suggestions. Also many thanks to Selin Ta\c skent, Yongsheng Zhang and Mustafa Kalafat for useful discussions and comments.

\end{document}